\documentclass[a4paper]{amsart}

\usepackage[latin1]{inputenc}
\usepackage{hyperref}
\usepackage{graphicx}
\usepackage{float}
\restylefloat{figure}

\newtheorem{theorem}{Theorem}
\newtheorem{corollary}{Corollary}
\newtheorem{proposition}{Proposition}
\newtheorem{lemma}{Lemma}
\newtheorem*{theorem*}{Theorem}

\newcommand{\remark}{\par\medskip\noindent\textbf{Remark.}\ }

\renewcommand{\epsilon}{\varepsilon}

\renewcommand{\Im}{\on{Im}}
\renewcommand{\Re}{\on{Re}}

\newcommand{\C}{\mathbb{C}}
\newcommand{\D}{\mathbb{D}}
\newcommand{\N}{\mathbb{N}}
\newcommand{\R}{\mathbb{R}}

\newcommand{\Z}{\mathbb{Z}}

\newcommand{\tend}{\longrightarrow}
\newcommand{\cst}{\on{cst}}

\newcommand{\setof}[2]{\big\{#1\,\big|\,#2\big\}}

\newcommand{\on}{\operatorname}
\newcommand{\cal}{\mathcal}

\newcommand{\wh}{\widehat}
\newcommand{\wt}{\widetilde}
\newcommand{\ds}{\displaystyle}

\begin{document}

\title{Semi-continuity of Siegel disks under parabolic implosion}
\author{Arnaud Chéritat}

\begin{abstract}We transfer the lower semi-continuity of Siegel disks with fixed Brjuno type rotation numbers to geometric limits. Here, we restrict to Lavaurs maps associated to quadratic polynomials with a parabolic fixed point.\end{abstract}

\maketitle

\tableofcontents

\section{Introduction}

We will prove here what was called hypothesis 4 in \cite{C}. There is no breakthrough here. The essential difficulty resides in getting a correct definition of perturbed horn maps, which we need to tend to the extended horn map of the limit parabolic point.

\section{Statement}

Let
\[P_\theta(z) = e^{i2\pi \theta} z + z^2.\]
The map $P_{p/q}$ with $p/q$ irreducible and $q>0$ has a parabolic point at $z=0$. The theory of parabolic implosion associates to this two families of (related) maps: the horn maps $h$, living in the cylinder and the Lavaurs maps $g$, living in the dynamical plane of $P_{p/q}$. Both depend on a parameter $\sigma\in\C$ called the \emph{phase}. Let us choose the phase such that the virtual multiplier of $h$ at one end of the cylinder is equal to $e^{i2\pi\theta}$, with $\theta$ a Brjuno number:
\[\theta\in\cal{B}.\]
To determine which end of the cylinder we will choose, we need to introduce a few definitions. 
Let us expand the $p/q$ in continued fraction:
\[\frac{p}{q} = a_0+\cfrac{1}{a_1+\cfrac{1}{\ddots+\cfrac{1}{a_m}}}\]
with $a_i\in\N^*$, which we more concisely write $p/q =[a_0,a_1,\ldots,a_m]$.
Let $\theta\in\cal B$ and
\[\theta_n = a_0+\cfrac{1}{\ddots+\cfrac{1}{a_m+\cfrac{1}{n+\theta}}} \]
Let us recall that $p/q$ has two such expansions, one of the form $[a_0,\ldots,a_k]$ with $a_k\geq 2$ and one of the form $[a_0,\ldots,a_k-1,1]$. In one case, $m=k$, in the other case $m=k+1$. We have
\[\theta_n \tend \frac{p}{q},\]
\[\on{sign}\big(\theta_n-\frac{p}{q}\big) = (-1)^m.\]
The end of the cylinder we consider in the upper end if $\theta_n>p/q$ and the lower end of $\theta_n<p/q$.

\begin{figure}[H]%
\begin{picture}(340,320)%
 \put(0,170){\scalebox{0.9}{\includegraphics{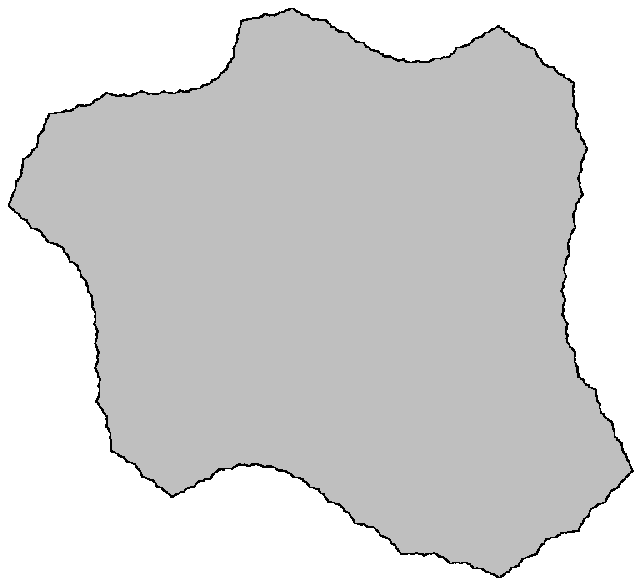}}}%
 \put(175,170){\scalebox{0.9}{\includegraphics{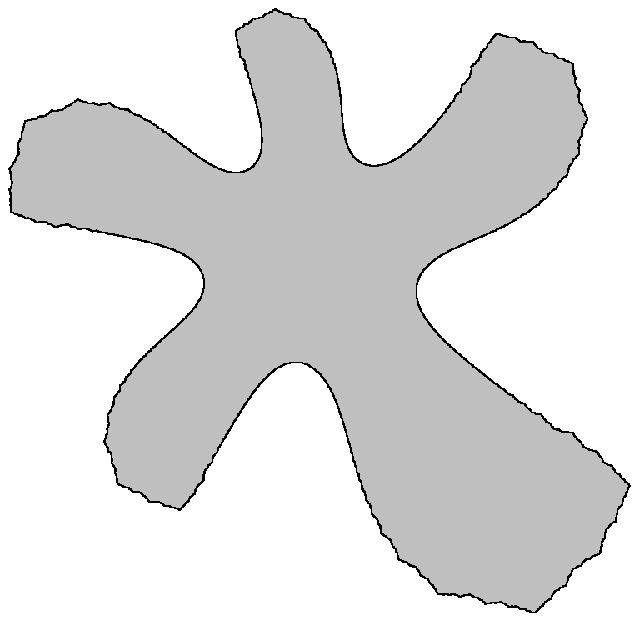}}}%
 \put(0,0){\scalebox{0.9}{\includegraphics{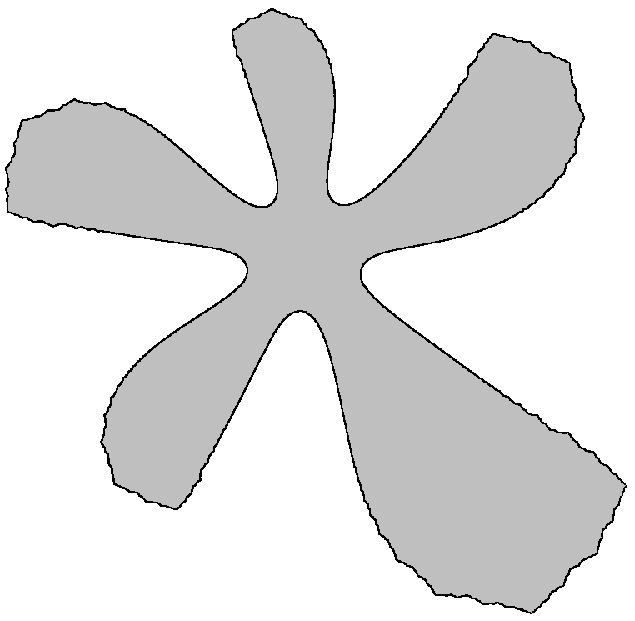}}}%
 \put(180,0){\includegraphics{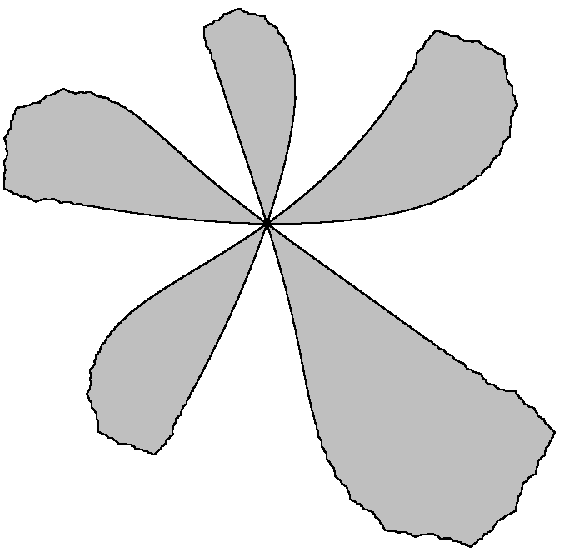}}%
\end{picture}%
\caption{Let $\theta_n=[0,2,2,n+\theta]$ with $\theta = \frac{\sqrt{5}-1}{2}$, the first three images show the Siegel disk of $P_{\theta_n}$ for $n=10,500,10000$, and the last is the virtual Siegel disk of $P_{2/5}$ they tend to. Here $\theta_n>p/q$.}%
\end{figure}

\begin{theorem*}
For all compact subset $C$ of the virtual Siegel disk, $\exists N$ such that $\forall n\geq N$, 
$C$ is contained in the Siegel disk of $P_{\theta_n}$.
\end{theorem*}

\remark The theorem extends to a more general setting than quadratic polynomials. 

\section{Tools of the proof}

\subsection{Semi-continuity of Siegel disks}

The following theorem is due to E.\ Risler, it is a corollary of \cite{R}. I wrote a specific proof in \cite{C}. Let $S_\theta$ be the set of analytic maps $f : \D \to \C$ fixing $0$ with multiplier $e^{i2\pi\theta}$ (we do not assume injectivity). The Siegel disk $\Delta(f)$ is the maximal subset of $\D$ containing $0$ and on which $f$ is conjugate to a rotation. Let
\[R_\theta(z) = e^{i2\pi \theta}z.\]

\begin{theorem}\label{thm_semi}
  $\forall \theta\in\cal{B}$, $\forall \epsilon>0$, $\exists \eta>0$ such that if $f\in S_\theta$ and $\|f-R_\theta\|_{\text{sup}}<\eta$ then $\Delta(f) \supset D(0,1-\epsilon)$.
  (For all Brjuno number, all maps in $S_\alpha$ sufficiently close to the rotation have a Siegel disk containing any given compact subset of $\D$.)
\end{theorem}

\subsection{Parabolic implosion}
\footnote{The theory of parabolic implosion (perturbation of parabolic points leading to geometric limits) has been developed by several authors, among which (alphabetical order)  A.~Douady, A.~Epstein, R.~Oudkerk, P.~Lavaurs, M.~Shishikura, and is still under development. Here we simply state the results that will be useful for our proof, after introducing the necessary definitions. We do not pretend to give an introduction to the topic.}

For all $n$ big enough, there are $q$ fixed points $z_1, \ldots, z_q$ of $P_{\theta_n}^q$ that tend to $0$ as $n\tend +\infty$. We have
\[z_i^q \underset{n\to +\infty}\sim \Big(\theta_n-\frac{p}{q}\Big)c\]
for some constant $c\in\C^*$, and moreover there is exactly one point $z_i$ close to each of the $q$-th roots of $(\theta_n-\frac{p}{q})c$. Thus they tend to $0$ asymptotically to $q$ equally spaced axes. We can label the axes with the symbols $1$ to $q$ and label the $q$ fixed points according to the axis they are close to. We will focus our attention on one of them, let us say $z_1$.
Let \[p(z) = z/|z| = e^{i \arg (z)}.\]

\begin{figure}[H]%
\begin{picture}(320,160)%
\put(0,20){%
 \put(0,0){\scalebox{0.3}{\includegraphics{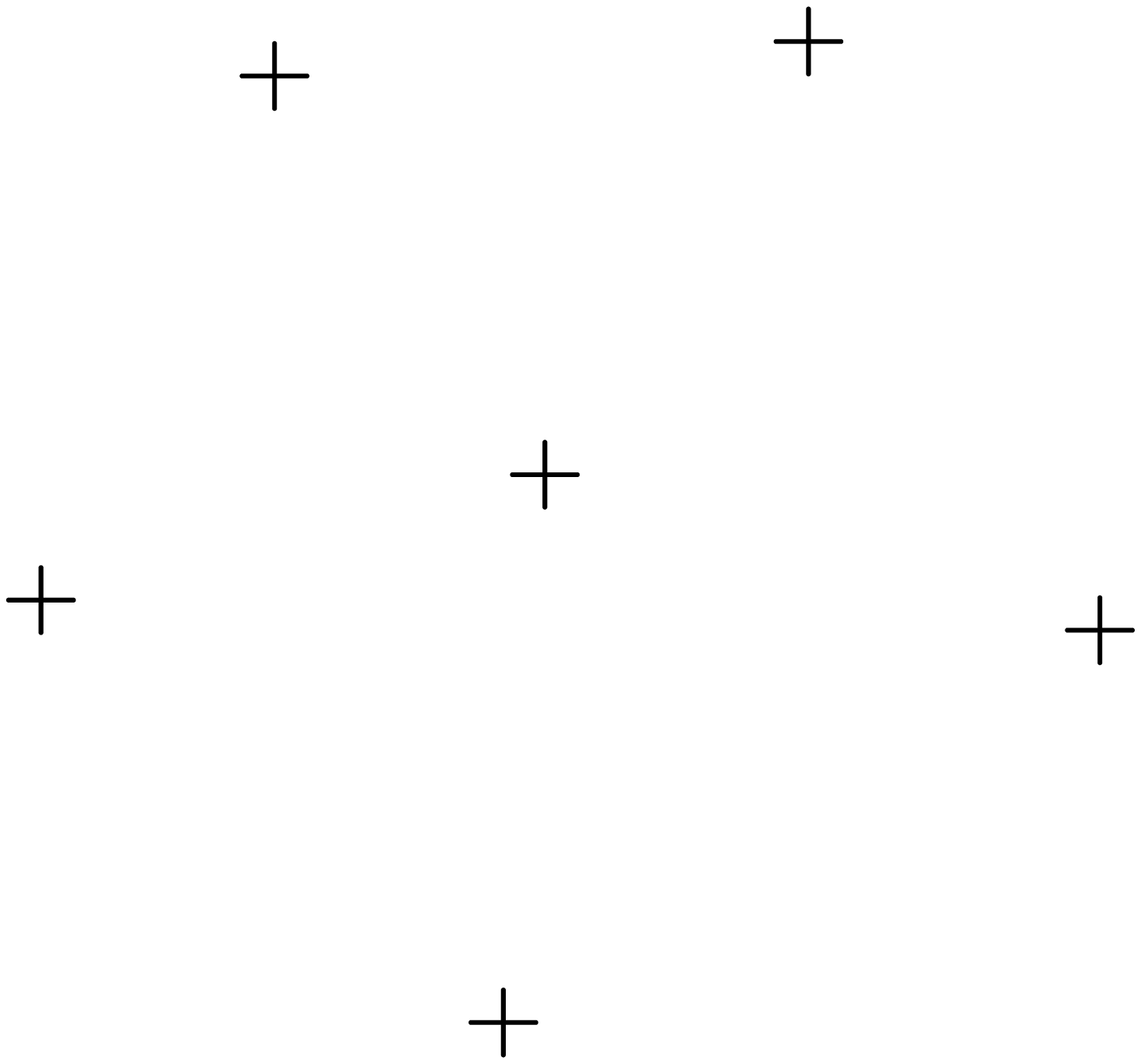}}}%
 \put(72,75){$0$}%
 \put(40,120){$z_1$}%
 \put(15,63){$z_2$}%
 \put(65,15){$z_3$}%
 \put(100,125){$z_4$}%
 \put(132,60){$z_5$}%
}%
\put(180,5){%
 \put(0,0){\scalebox{0.3}{\includegraphics{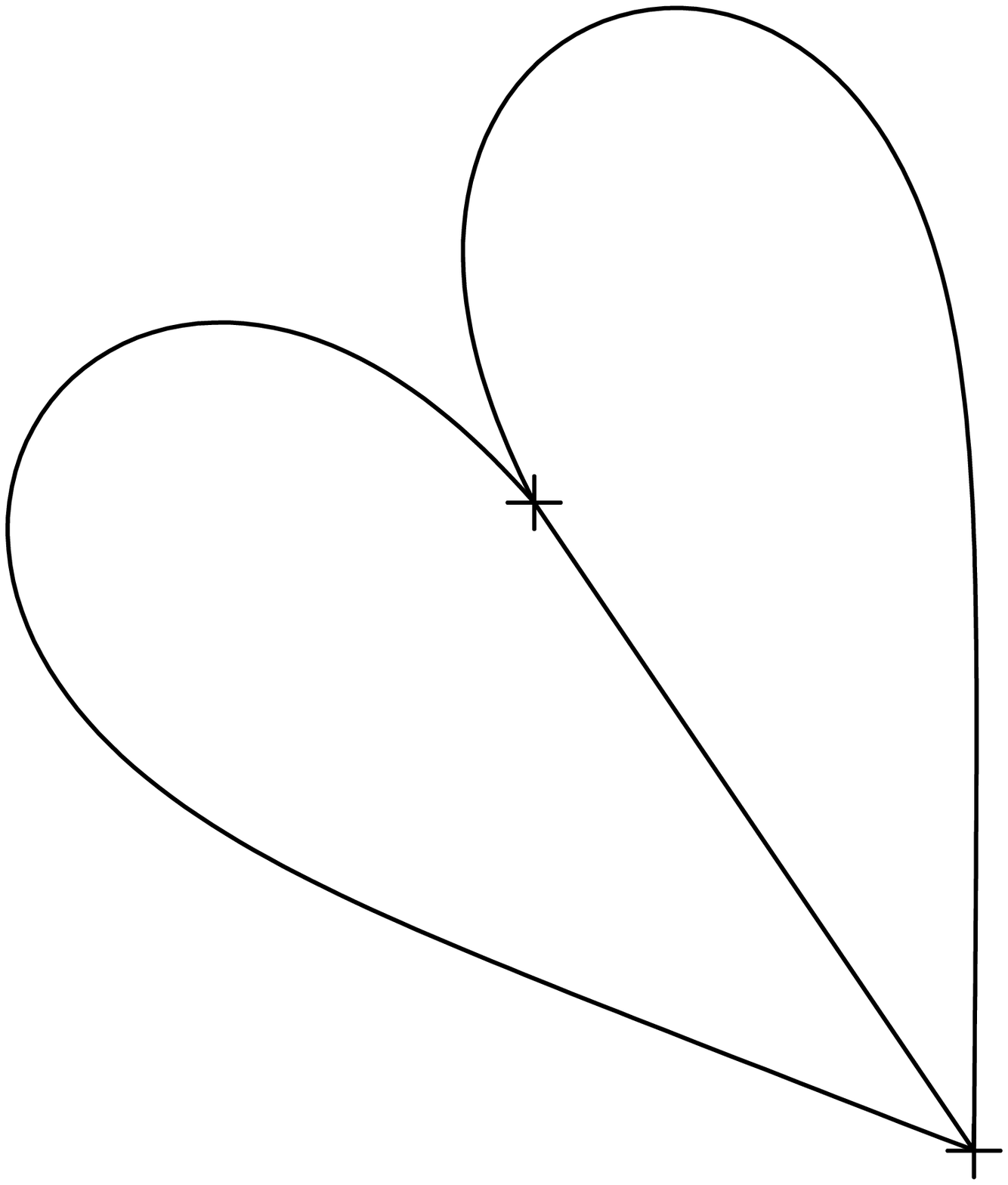}}}%
 \put(135,0){$0$}%
 \put(77,86){$z_1$}%
 \put(90,110){$P_n^-$}%
 \put(30,70){$P_n^+$}%
 \put(101,48){$I_n$}%
}%
\end{picture}%
\caption{The fixed points and the petal $P_n$.}%
\end{figure}%

We will work on a small ball $B(0,r)$, where $r$ will be fixed later. Let us consider the map $\on{pow}_n : z\mapsto (z/p(z_1))^q$ defined on the sector $``\arg(z/z_i) \in ]-\frac{\pi}{q},\frac{\pi}{q}["$ (so that $\on{pow}_n$ is injective). Let us consider the two biggest open disks whose boundary passes by $0,1$ and that are contained in the image of $B(0,r)\cap S$ by $\on{pow}_n$, i.e. in $B(0,r^q)\setminus]-\infty,0]$. Let $P_n$ be the preimage by $\on{pow}_n$ of the union of these two disks (see the picture). This is the petal associated to $z_1$. Let us cut the union of two disks by the real line. It corresponds to a cut $I_n = [0,z_1]$ of $P_n$. If $\theta_n>p/q$, we will call $P_n^+$ the preimage of the part above the real line, and $P_n^-$ the other one. If $\theta_n<p/q$ this will be the other way round.

Let $\zeta_n : P_n \to \C$ be defined by
\[\zeta_n(z) = \frac{-1}{\scriptstyle i2\pi q^2(\theta_n-p/q)}\log\Big(1-(z_1/z)^q\Big).\]
The map $1-z\mapsto(z_1/z)^q$ sends $P_n$ to a sector centered on $0$, on which a branch of $\log$ is chosen.
The map $\zeta_n$ is univalent and sends $P_n$ to a vertical band, mapping $0$ to the upper end and $z_1$ to the lower end if $\theta_n>p/q$ and the other way round otherwise.

\begin{figure}[H]%
\begin{picture}(340,230)%
\put(30,0){\scalebox{0.8}{%
 \put(-70,30){%
  \put(-10,0){\scalebox{0.4}{\includegraphics{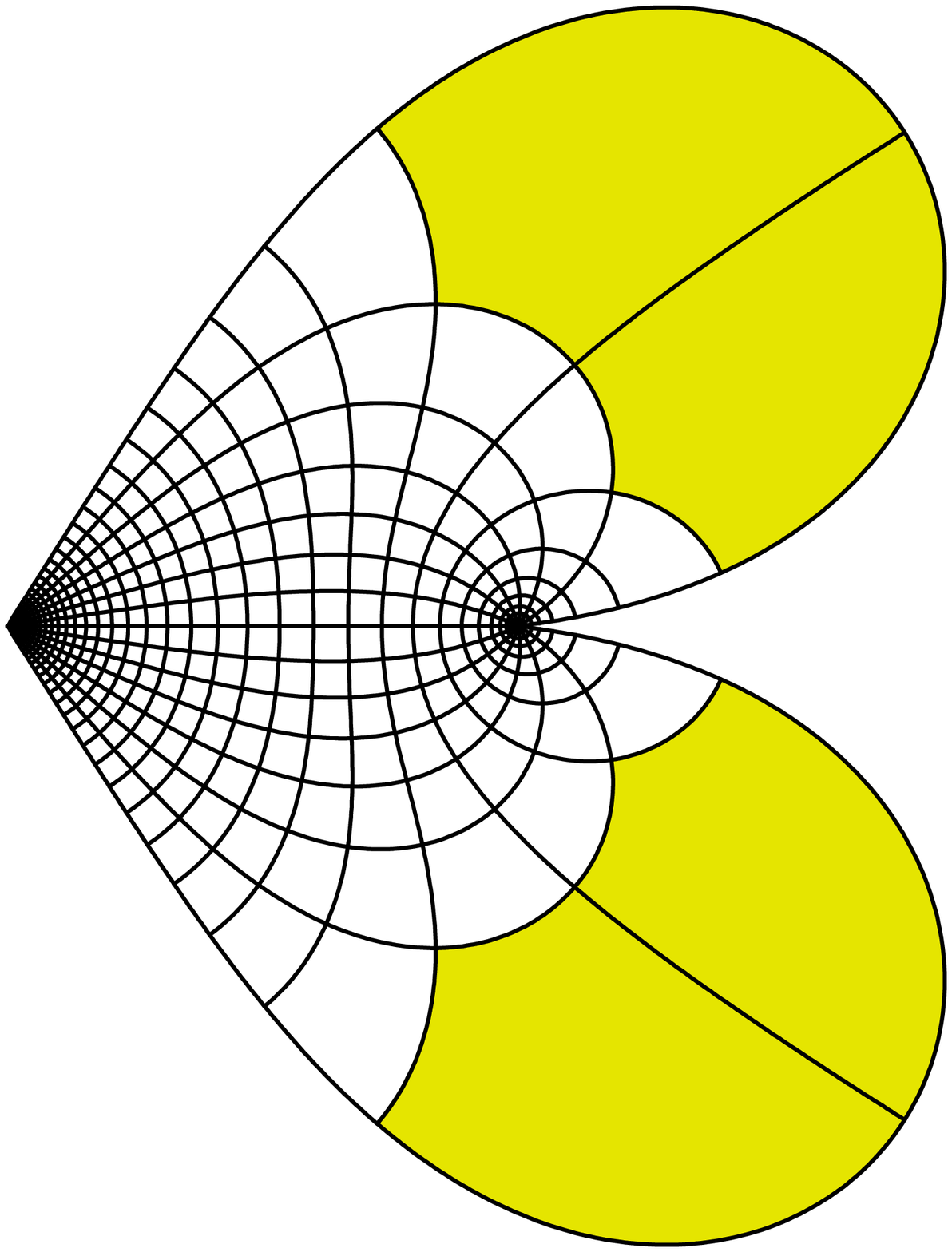}}}%
  \put(170,28){\scalebox{0.4}{\includegraphics{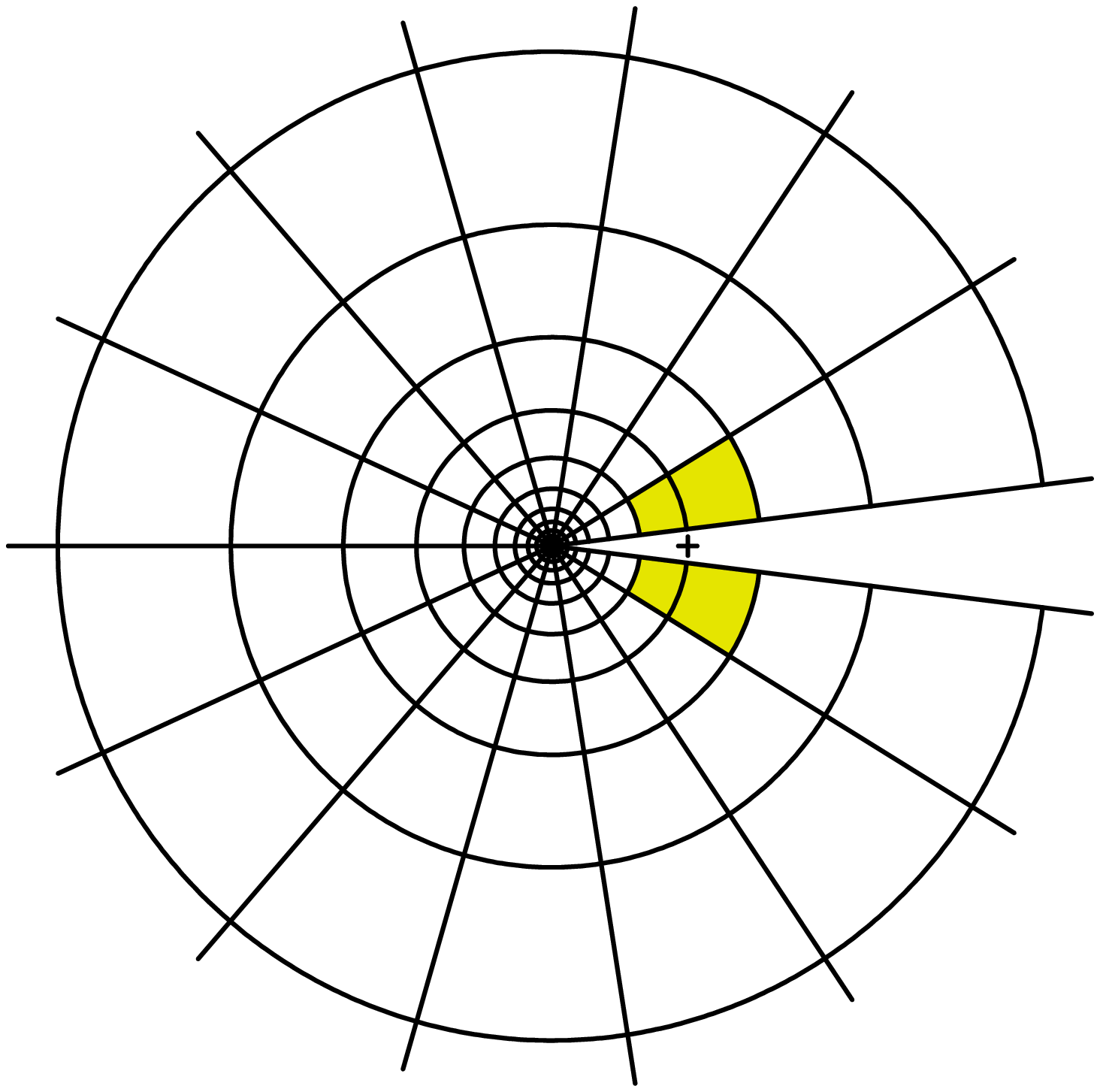}}}%
 }%
 \put(285,0){\scalebox{0.5}{\includegraphics{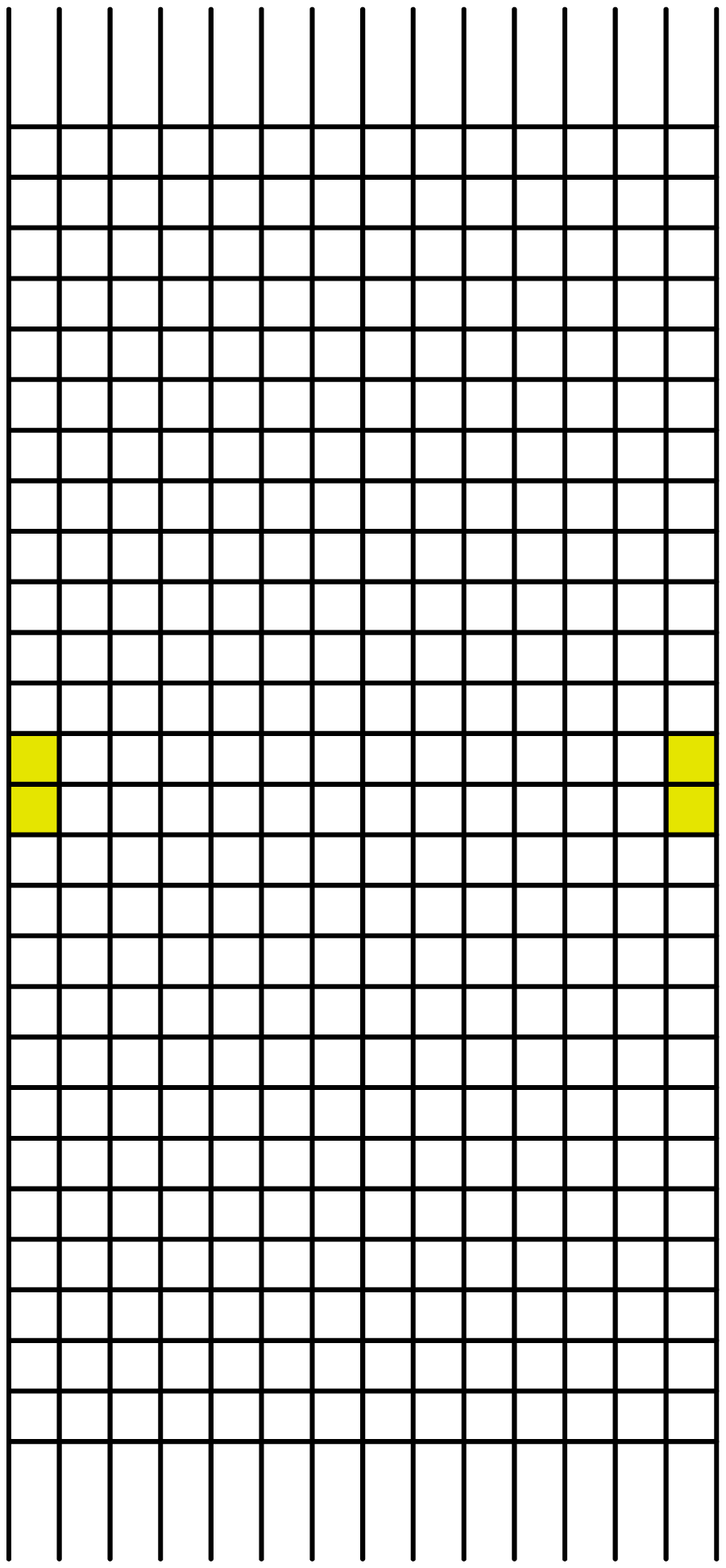}}}%
}}%
\end{picture}%
\caption{The petal (rotated), its image under the map $z\mapsto 1-\frac{z_1}{z^q}$, and its image under $\zeta_n$. The grid and the yellow parts correspond to each other. In the middle image, the yellow boxes are close to the point $z=1$ which is marked by a cross.}%
\end{figure}

Let $\alpha_1$ be the argument of the axis $z_1$ is close to. Let $\on{pow} : z \mapsto (z/e^{i\alpha_1})^{q}$ be defined on the sector $``\arg(z/e^{i\alpha_1}) \in ]-\frac{\pi}{q},\frac{\pi}{q}["$. Let $P$ be the preimage by $\on{pow}$ of the union of the two biggest disks tangent to $\R$ at $0$ and contained in $B(0,r^q)$ (these are the limits of the previous disks). Then every compact subset of $P$ is eventually contained in $P_n$. Let $P_+$ and $P_-$ be the analogs of $P_n^+$ and $P_n^-$. Then, provided $r$ is small enough, $P_+$ is a repelling petal and $P_-$ is an attracting petal for $P_{p/q}$.

\begin{lemma}\label{lem_zetalim}
  There exist $k_+,k_-\in\Z$ (with $k_- - k_+ = (-1)^m$) such that
  \[\zeta_n(z) + \frac{k_+}{q^2(\theta_n-p/q)} \tend \frac{c}{i2\pi q^2 z^q}\]
  on every compact subset of $P^+$ and
  \[\zeta_n(z) + \frac{k_-}{q^2(\theta_n-p/q)} \tend \frac{c}{i2\pi q^2 z^q}\]
  on every compact subset of $P^-$.
\end{lemma}

\begin{figure}[H]%
\begin{picture}(130,180)%
\put(0,0){\scalebox{0.3}{\includegraphics{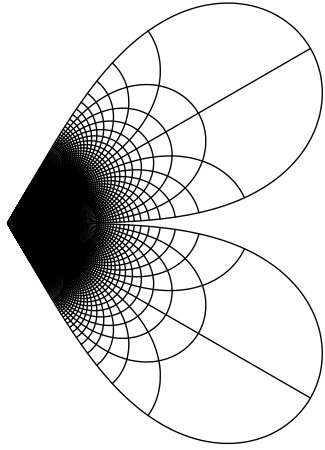}}}%
\end{picture}%
\caption{The limit petal $P_+ \cup P_-$ (rotated), decorated by the preimage of the orthogonal grid by the limit map $z\mapsto \frac{c}{i2\pi q^2 z^q}$.}%
\end{figure}%

Let $\pi : \C \to \C/\Z$ be the quotient map. Then $\exists r_0>0$, $\forall r<r_0$, $\exists N\in\N$, $\forall n\geq N$, the following two theorems hold.

\begin{figure}[H]%
\begin{picture}(215,185)%
\put(0,0){\scalebox{0.5}{\includegraphics{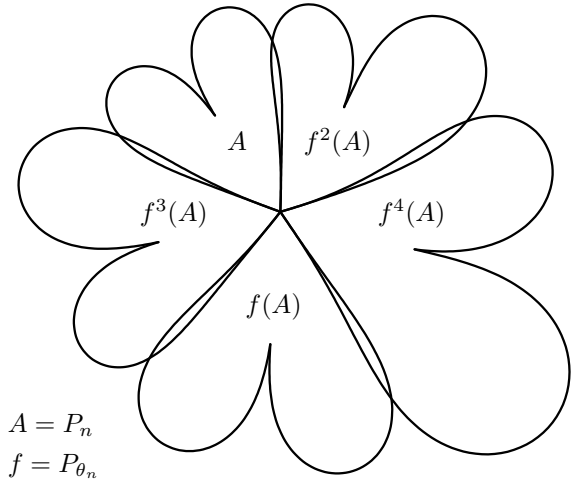}}}%
\put(83,127){$A$}%
\put(90,65){$f(A)$}%
\put(112,127){$f^2(A)$}%
\put(50,100){$f^3(A)$}%
\put(140,100){$f^4(A)$}%
\put(0,20){$A=P_n$}%
\put(0,5){$f=P_{\theta_n}$}%
\end{picture}%
\caption{This picture illustrates the domain $P_n$ and its first $4$ iterates by $P_{\theta_n}$. Notice the overlaps.}%
\end{figure}

\newcounter{savedCounter}
\begin{theorem}\footnote{Points (\ref{item_disj_m}) through (\ref{item_trap}) are elementary consequence of $P_{\theta_n}$ being close to $R_{\theta_n}$ near $z=0$. Points (\ref{item_entry}) through (\ref{item_overlap}) are consequences of the following estimate: $\Big|\zeta_n\circ P_{\theta_n}^q(z) \circ\zeta_n^{-1}-(z+1)\Big| \leq \epsilon$ for $r$ small enough (and $N$ big enough, depending on $r$). Point (\ref{item_ptom}) is a consequence of the same estimate for an other branch of $\zeta_n$.}
\begin{enumerate}
\item\label{item_disj_m} The following sets are pairwise disjoint: $P_n^-$, $P_{\theta_n}(P_n^-)$, \ldots, $P_{\theta_n}^{q-1}(P_n^-)$.
\item\label{item_disj_c} The closures of two of them intersect only at $0$.
\item\label{item_disj_lm} The sames hold for $P_-$, \ldots, $P_{p/q}^{q-1}(P_-)$.
\item\label{item_disj_i} If $z\in P_n^+$ and $0<k<q$ then $P_{\theta_n}^k(z) \notin I_n$. 
\item\label{item_trap} $\exists r'<r$ ($r'$ depends on $n$) such that every point in $B(0,r')$ eventually falls in $P_n^-$ under iteration of $P_{\theta_n}$.
\item\label{item_entry} For all $z\in\partial P_n^-\setminus I_n$, $P_{\theta_n}^{q}(z) \in P_n^-$.
\item\label{item_mtop} Every point of $P_n^-$ is eventually mapped to $P_n^+$ under iteration of $P_{\theta_n}^q$. 
\item\label{item_overlap} If $z\in P_n^- \cap P_{\theta_n}^k(P_n^+)$ with $0\leq k<q$ then $P_{\theta_n}^q(z) \in P_n^-$.
\item\label{item_ptom} Note $]a,b[\times\R = \zeta_n(P_n)$. There exists $M>0$ (which depends on $r$ but not on $n$) such that every point $z \in P_n^+$ with $\Im(\zeta_n(z))>M$ and $\Re(\zeta_n(z))>b-2$ is eventually mapped to $P_n^-$ under iteration of $P_{\theta_n}$. 
\setcounter{savedCounter}{\value{enumi}}
\end{enumerate}
\end{theorem}

\begin{theorem}[perturbed Fatou coordinates]\ % LaTeX does not newline without a char here
\begin{enumerate}
\setcounter{enumi}{\value{savedCounter}}
\item\label{item_phin} There exists a univalent function $\Phi_n : P_n \to \C$ such that $\Phi_n \circ P_{\theta_n}^q = T_1 \circ \Phi_n$ holds on $P_n \cap (P_{\theta_n}^q)^{-1}(P_n)$.
\item\label{item_deriv} The derivative of $\Phi_n \circ \zeta_n^{-1}$ stays in $B(1,1/4)$.
\item\label{item_surj} The image of $P_n^+$ under $\Phi_n$ is surjectively mapped to $\C/\Z$ by $\pi$. The same holds for $P_n^-$.
\item\label{item_cvx} If both $z$ and $z+n$ belong to $\Phi_n(P_n^+)$, then $z+k \in \Phi_n(P_n^+)$ for $0 < k < n$. The same holds for $P_n^-$.
\item\label{item_nu} $\ds \on{sign}(\theta_n-p/q).\on{Im}(\Phi_n(z)) \tend +\infty \iff z\tend 0$ 
\item\label{item_dev} $\ds \Phi_n(z) \underset{z\tend 0}{=} \frac{\log z}{i2\pi(q\theta_n-p)} + \cst + o(1)$.
\item\label{item_lim} For any given $a\in P_+$ (resp.\ $P_-$), $\Phi_n(z)-\Phi_n(a)$ tends uniformly on every compact subset of $P_+$ (resp.\ $P_-$) to a Fatou coordinate for $P_{p/q}$. 
\end{enumerate}
\end{theorem}

Points (\ref{item_surj}), (\ref{item_cvx}) and (\ref{item_nu}) follow from (\ref{item_deriv}).

Point (\ref{item_lim}) motivates the following definition: we fix any $a_+ \in P_+$ and $a_- \in P_-$, and define $\Phi_+ : P_+ \to \C$ to be the limit of $\Phi_n(z)-\Phi_n(a_+)$ for $z\in P_+$ and analogously for $\Phi_- : P_- \to \C$.

\section{Proof}

\subsection{Defining perturbed horn maps}

This definition requires some care.

We first define a function $k$ on $\C$:
\[k(z) = \inf\setof{k\in\N}{P_{\theta_n}^k(z) \in P_n^-} \in \N\cup\{+\infty\}.\]
Since $P_n^-$ is open, $k$ is upper semi-continuous.

Then we let
\[U=\setof{z\in P_n^+}{k(z)<+\infty \text{ and }
P_{\theta_n}^k(z) \notin I_n
\text{ for } 0\leq k \leq k(z)}.\]
The set $U$ is open and $U \subset P_+$.

\medskip

From (\ref{item_disj_i}) and (\ref{item_overlap}), we get:
\begin{lemma}\label{lem_U}
  If $z\in P_n^+$ and $P_{\theta_n}^q(z) \in P_n^+$ then $z\in U \iff P_{\theta_n}^q(z) \in U$.
\end{lemma}

This, together with~(\ref{item_phin}) and~(\ref{item_cvx}) implies that $\Phi_n(U) = \Phi_n(P_n^+) \cap \pi^{-1}(\wt{U})$ for some subset $\wt{U}$ of $\C/\Z$.

We set
\[\on{def}(h_n) = \wt{U}\]
which is an open subset of $\C/\Z$.
For $w\in \wt{U}$ we let $w'\in \Phi_n(P_n^+)$ be any representative of $w$.
Then $w'=\Phi_n(z)$ for a unique $z\in U$ and we let
\[h_n(w) = \Phi_n(P_{\theta_n}^{k(z)}(z)) \bmod \Z \in \C/\Z.\]

\noindent\textbf{Claim:} $h_n(w)$ is independent of the representative $w'$ of $w$ we chose.
\proof Indeed, if $w_1=\Phi_n(z_1)$ and $w_2=w_1+m=\Phi_n(z_2)$ are two representatives, then
(\ref{item_cvx}) and (\ref{item_phin}) imply $P_{\theta_n}^{jq}(z_1) \in P_n^+$ for $0\leq j\leq m$ and $z_2=P_{\theta_n}^{mq}(z_1)$. Let $u_k=P_{\theta_n}^k(z_1)$.
\par\noindent Case 1: $k(z_1)\geq mq$. Then $k(z_2)=k(z_1)-mq$ and thus
$P_{\theta_n}^{k(z_2)}(z_2) = P_{\theta_n}^{k(z_1)}(z_1)$.
\par\noindent Case 2: $k(z_1)\leq mq$. Then point (\ref{item_overlap}) implies $u_{k(z_1)+jq} \in P_n^-$ for all $j$ such that $k(z_1)\leq k(z_1)+jq<mq+q$. By (\ref{item_disj_m}), these are the only values of $k$ such that $k(z_1)\leq k <mq+q$ and $u_k \in P_n^-$. In particular,
$k(z_2)=b$ where $b\in[0,q[$ is the remainder of the Euclidean division $k(z_1)=aq+b$.
By point (\ref{item_phin}), we have $\Phi_n(u_{k(z_1)+jq}) = j + \Phi_n(u_{k(z_1)})$.
Whence $\Phi_n(P_{\theta_n}^{k(z_2)}(z_2)) = \Phi_n(u_{mq+b}) = m-a+\Phi_n(u_{k(z_1)})
=m-a+\Phi_n(P_{\theta_n}^{k(z_1)}(z_1))$.
\qed

\par\medskip
\noindent\textbf{Claim.} The map $h_n$ is analytic. 
\proof It is obvious for points $w\in\wt{U}$ such that $k$ is locally constant at $z$. If $w_0 = \pi(\Phi_n(z_0))$ with $z_0\in U$ is such that $k$ is not locally constant at $z_0$, it means that there is $k'<k(z_0)$ such that $P_{\theta_n}^{k'}(z_0)\in\partial P_n^-$. Since $I_n$ must be avoided, point (\ref{item_entry}) says $P_{\theta_n}^{k'+q}(z_0) \in P_n^-$. Now by (\ref{item_disj_m}), $k(z_0)=k'+q$. In particular, $k'$ is unique. By (\ref{item_disj_c}), one concludes that in a neighborhood $V$ of $z_0$, $k(z) \in \{k', k'+q\}$. We can choose $V$ small enough so that $P_{\theta_n}^{k'+q}(V)\subset P_n^-$. For $z\in V$ so that $k(z)=k'$, (\ref{item_phin}) implies $\Phi_n(P_{\theta_n}^{k'+q}(z))=\Phi_n(P_{\theta_n}^{k'}(z)) \bmod \Z$. Thus $h_n(w) = \Phi_n(P_{\theta_n}^{k'+q}(z)) \bmod \Z$ holds on a neighborhood of $w_0$. \qed

\medskip

By (\ref{item_mtop}), we get:
\begin{proposition}
  $\forall z\in P_n$, if $\Phi_n(z) \in \on{def}(h_n) \bmod \Z$ then $\exists k \in \N$ such that\footnote{$k>0$ is an important point} $k>0$, $z'=P_{\theta_n}^k(z) \in P_n^+$ and $\Phi_n(z') = h_n(\Phi_n(z)) \bmod \Z$.
\end{proposition}

\begin{corollary}\label{cor_inforb}
  $\forall z\in P_n$, if $\Phi_n(z) \bmod \Z$ has an infinite orbit under $\on{def}(h_n)$, then the orbit of $z$ under $P_{\theta_n}$ passes an infinite number of times in $P_n^+$.  
\end{corollary}

\medskip

By (\ref{item_surj}), (\ref{item_nu}) and (\ref{item_trap}), the set $\on{def}(h_n)$ is a neighborhood of $\nu$, and $h_n$ has a continuous (and thus analytic) extension $\wh{h}_n$ fixing it.\footnote{We consider the cylinder completed by its two ends as a Riemann surface isomorphic to the Riemann sphere (via the map $(z\bmod \Z)\mapsto e^{i2\pi z}$).}

Let $\nu$ be the upper end if $\theta_n>p/q$ and the lower end if $\theta_n<p/q$.

\begin{proposition}[folk.]\label{prop_hmul}
  The multiplier of $\wh{h}_n$ at the end $\nu$ is equal to
  \[\exp\big((-1)^{m+1}i2\pi\theta\big).\]
\end{proposition}
\noindent This is a well known corollary of point (\ref{item_dev}), of a local combinatorial study, and of the theory of continued fractions.

\bigskip

The (unperturbed) horn map $h$ can be defined as follows: given $w\in\C/\Z$ let $w'\in \Phi_+(P_+)$ be any representative of $w$ in $\C$ and $z=\Phi_+^{-1}(w')$. If there is a $k \in \N$ such that $z'=P_{\theta_n}^k(z) \in P_-$, then we define 
\[h(w)=\Phi_-(z') \bmod \Z.\]
Otherwise we decide that $z \notin \on{def}(h)$.
Analogously, the map $h$ is well defined, the set $\on{def}(h)$ is open, $h$ is analytic, $\on{def}(h)$ is a neighborhood of $\nu$, and $h$ has an analytic continuation $\wh{h}$ fixing it (cf.~\cite{DH}, \cite{E}, \cite{S}, etc\ldots). 

\begin{lemma}\label{lem_cvh}
  \[ T_{-\Phi_n(a_-)} \circ \wh{h}_n \circ T_{\Phi_n(a_+)} \underset{n\to+\infty}{\tend} \wh{h}\]
  in the sense described above.
\end{lemma}

\proof Let $w\in\on{def}(\wh{h})$.

\smallskip

\noindent Case 1: $w\neq \nu$. By definition and (\ref{item_lim}), $T_{-\Phi_n(a_+)} \circ \Phi_n \tend \Phi_+$ uniformly on every compact subsets of $P_+$. Thus $\Phi_n^{-1}\circ T_{\Phi_n(a_+)} \tend \Phi_+^{-1}$ uniformly on every compact subsets of $\Phi_+(P_+)$. If $z\in P_+$ is eventually mapped to $P_-$ under iteration of $P_{p/q}$, say $P_{p/q}^k(z) \in P_-$, then there is a neighborhood $V$ of $z$ and $N$ such that $\forall n\geq N$, $P_{\theta_n}^k(V) \subset P_n^-$. By taking a slightly smaller $V$, we may assume also that
$0\notin\overline{P_{p/q}^k(V)}$. Since $I_n \tend \{0\}$ as $n \tend +\infty$, we know that for all $n$ big enough, $V$ avoids $I_n$ before getting in $P_n^-$ under $P_{\theta_n}$ (otherwise, $\overline{P_{p/q}^k(V)}$ would contain $0$). Therefore, there is a neighborhood of $w$ that is contained in $\on{def}(h_n)$ for $n$ big enough. Now, either $z$ does not fall in $\partial P_-$ before getting in $P_-$, in which case there is a neighborhood $V'\subset V$ of $z$ and $N'\geq N$ such that $\forall n\geq N'$, $\forall z'\in V'$, $k_n(z') = k(z)$. Or $\exists k'<k(z)$ such that $P_{p/q}^{k'}(z)\in\partial P_-$. Since $P_-$ is a petal, we know by (\ref{item_disj_lm}) that $k(z)=k'+q$, whence $k'$ is unique. This implies, by (\ref{item_disj_lm}) and (\ref{item_phin}) that we have $\wh{h}_n(w) = \Phi_n(P_{\theta_n}^{k'+q}(z)) \bmod \Z$ which holds for all $n$ big enough on a neighborhood of $w$ independent of $n$. In both cases, the convergence then follows from the definitions and from $T_{-\Phi_n(a_-)}\circ\Phi_n \tend \Phi_-$. 

\medskip

\noindent Case 2: $w=\nu$. We will prove that there is a neighborhood $V$ of $\nu$ that is eventually contained in $\on{def}(\wh{h}_n)$. Then, using case 1 on a circle $A\subset V$ around $\nu$, we deduce that the limit $h_n \tend h$ must be uniform on $A$. But since $\wh{h}_n(\nu)=\nu$, the limit $\wh{h}_n \tend \wh{h}$ must be uniform above $A$ for the spherical metrics, which ends the proof of case 2. For the existence of $V$, consider the set $U$ of points in the band $\zeta_n(P_n) = ]a,b[ \times \R$ whose real part are $\geq b-2$ and imaginary part are $\geq M$ with $M$ given by (\ref{item_ptom}). Since $\zeta_n$ converges on compact subsets of $P_+$ (cf.\ lemma~\ref{lem_zetalim}), it converges at $a_+$, and this implies by (\ref{item_deriv}) that the sets $\pi\circ T_{-\Phi_n(a_+)}\circ\Phi_n\circ \zeta_n^{-1}(U)$ contain a common\footnote{i.e.\ independent of $n$} neighborhood of $\nu$ in $\C/\Z$. Let us call it $V$. By (\ref{item_ptom}), $V \subset \on{def}(\wh{h}_n)$.

\qed

Since $\wh{h}$ has a non zero multiplier at the end $\nu$,
we have $\wh{h}(z)-z \tend \tau$ as $z\tend\nu$. Together with proposition~\ref{prop_hmul} this gives:
\[\Phi_n(a_+)-\theta-\Phi_n(a_-) \underset{n\to+\infty}{\tend} \tau\]
Thus we can rephrase lemma~\ref{lem_cvh} in a more convenient way:
\[ T_{-\Phi_n(a_+)} \circ \wh{h}_n \circ T_{\Phi_n(a_+)} \underset{n\to+\infty}{\tend} 
T_{-\theta-\tau}\circ \wh{h}.\]
In other words
\[ \wh{f}_n \underset{n\to+\infty}{\tend} \wh{f}\]
with
\[ \wh{f}_n = T_{-\Phi_n(a_+)} \circ \wh{h}_n \circ T_{\Phi_n(a_+)}\]
and
\[ \wh{f} = T_{-\theta-\tau}\circ \wh{h},\]
and $\wh{f}$ and $\wh{f}_n$ have the same multiplier at $\infty$ (the one given by proposition~\ref{prop_hmul}). Since this number is a Brjuno number, we can apply the semi-continuity result (theorem~\ref{thm_semi}) and deduce from this: 

\begin{lemma}\label{lem_sch}
Every compact subset of $\Delta(\wh{f}_n)$ is eventually contained in $\Delta(\wh{f})$. Note that we include the end $\nu$ in these Siegel disks, so in particular, there is a neighborhood of $\nu$ that is eventually contained in $\Delta(\wh{f}_n)$.
\end{lemma}

\subsection{Transferring the semi-continuity}

Let us define the shorthands $\Delta(f) = \pi^{-1}(\Delta(\wh{f}))$, and $\Delta(f_n) = \pi^{-1}(\Delta(\wh{f}_n))$.
The \emph{virtual Siegel disk} is defined\footnote{There are several equivalent definitions of the virtual Siegel disk. For instance we may have used the attracting Fatou coordinate $\Phi_-$ instead of $\Phi_+$. We may also have used the maximal extension of $\Phi_+^{-1}$ instead of taking forward images by $P_n$. Or the maximal extension of $\Phi_-$ together with covering properties of $\Phi_-$.} as
\[\Delta_\infty = \bigcup_{k\in\N} P_{p/q}^k (\Phi_+^{-1}(\Delta(f))).\]
Let $\Delta_n$ be the (usual) Siegel disk of $P_{\theta_n}$.\footnote{Linearizability of $P_{\theta_n}$ can be justified in several ways. Invariance of $\cal B$ under $\on{PSL}(\Z^2)$ and $\theta\in\cal{B}$ imply that $\theta_n \in \cal{B}$. Or one may argue that $\wh{h}_n$ is a (Yoccoz) renormalization (high enough in the cylinder) of $P_{\theta_n}$ and that linearizability of $\wh{h}_n$ thus implies that of $P_n$.}
By lemma~\ref{lem_sch}, $\exists y_0>0$, $\exists N\in\N$, such that $\forall n\geq N$, $\Delta(f_n)$ contains $``\Im(z)>y_0"$.

\medskip\noindent\textbf{Claim:} $\forall n\geq N$ and $\forall z\in P_n$, if $\Im(\Phi_n(z)-\Phi_n(a_+))> y_0$ then $z \in \Delta_n$.

\proof 
By (\ref{item_deriv}), $\Phi_n(P_n) \cap ``\Im(z)>\text{something}"$ is connected. Thus
$\Phi_n^{-1}(``\Im(z-\Phi_n(a_+))>y_0")$ is a connected open set whose points have by corollary~\ref{cor_inforb} an orbit under $P_{\theta_n}$ which passes an infinite number of times in the bounded set $P_n^+$. Therefore, it is contained in the Fatou set of the polynomial $P_{\theta_n}$. Using (\ref{item_deriv}) again, we see that this open set adheres to $0$. Since $P_{\theta_n}$ is linearizable at $0$, the connected open set is contained in the Siegel disk.
\qed

\medskip

Let $z_0$ be in the virtual Siegel disk. Let us write $z_0=P_{p/q}^{k_0}(\Phi_+^{-1}(w_0))$ with $w_0 \bmod \Z \in \Delta(f)$. Let $\gamma : [0,1] \to \Delta(f)$ a path from $w_0$ to any point $w_1$ in $``\Im(z)>y_0"$. By lemma~\ref{lem_sch}, the curve $\gamma$ is contained in $\Delta(f_n)$ for $n$ big enough. Since this curve is compact, there exists $m_0>0$ such that $T_{-m_0}\circ\gamma([0,1]) \subset \Phi_+(P_+)$. By definition and (\ref{item_lim}), there exists a neighborhood $U$ of $\gamma([0,1])$ such that $T_{-m_0}(U) \subset T_{-\Phi_n(a_+)}\circ\Phi_n(P_n)$ for $n$ big enough. Up to taking a smaller neighborhood, we will choose $U$ open, connected and contained in $\Delta(f)$. Then, by corollary~\ref{cor_inforb}, $\Phi_n^{-1}\circ T_{\Phi_n(a_+)} \circ T_{-m_0} (U)$ is contained in the Fatou set, and by the claim above, it contains a point in $\Delta_n$ (the image of $\gamma(1)$). Being moreover connected, it is contained in $\Delta_n$. Therefore \[P_{\theta_n}^{k_0} \circ \Phi_n^{-1}\circ T_{\Phi_n(a_+)} \circ T_{-m_0} (U) \subset \Delta_n.\]
Since the limit $P_{p/q}^{k_0} \circ \Phi_+^{-1} \circ T_{-m_0}$ of the functions $P_{\theta_n}^{k_0} \circ \Phi_n^{-1}\circ T_{\Phi_n(a_+)} \circ T_{-m_0}$ is non constant, the image of $U$ and contains some neighborhood of $z_0$ for $n$ big enough. This ends the proof of the main theorem.

\bibliographystyle{plain}

\end{document}